\documentclass{article}
\usepackage[utf8]{inputenc}
\usepackage[english]{babel}
\usepackage{amscd,amsmath,amssymb,indentfirst,a4wide}

\newcounter{Mycounter}[section]
\newcounter{lemma}[section]
\newcounter{claim}[section]
\newcounter{corollary}[section]
\newcounter{theorem}[section]
\renewcommand{\thetheorem}{{\thesection.\arabic{theorem}}}
\newcommand{\theorem}{%
     \setcounter{theorem}{\value{Mycounter}}
     \refstepcounter{theorem}
     \stepcounter{Mycounter}
     {\bf Theorem \thetheorem:\ }}

\newcounter{conjecture}[section]
\newcounter{proposition}[section]
\newcounter{definition}[section]
\renewcommand{\thedefinition}
       {{\thesection.\arabic{definition}}}
\newcommand{\definition}{
     \setcounter{definition}{\value{Mycounter}}
     \refstepcounter{definition}
     \stepcounter{Mycounter}
     {\bf Definition \thedefinition:\ }}

\newcounter{example}[section]
\newcounter{remark}[section]
\newcounter{problem}[section]
\newcounter{question}[section]
\makeatletter
\@addtoreset{equation}{section}
\@addtoreset{footnote}{section}
\makeatother

\def\tilde{\widetilde}

\newcommand{\im}{\operatorname{im}}

\begin{document}
\begin{center}
\begin{Large}
\bf{Curves on Oeljeklaus-Toma Manifolds}
\end{Large}\\[5mm]
\begin{large}
Sima Verbitsky
\end{large}\\[6mm]
\end{center}

{\small 
\hspace{0.15\linewidth}
\begin{minipage}[t]{0.6\linewidth}
\begin{center}
\bf Abstract
\end{center}
\vspace{-0.03\linewidth}
Oeljeklaus-Toma manifolds are complex non-K\"ahler manifolds constructed by Oeljeklaus and Toma from
certain number fields, and generalizing the Inoue surfaces $S_m$. We prove that Oeljeklaus-Toma manifolds contain no compact complex curves.
\end{minipage}
}

\tableofcontents

\section{Introduction}

Oeljeklaus-Toma manifolds (defined in \cite{_Oeljeklaus_Toma_}) are compact complex manifolds that are 
a generalization of Inoue surfaces
(defined in \cite{_Inoue_}). Let us describe them in detail.

\subsection{Oeljeklaus-Toma manifolds}

Let $K$ be a number field (i.e. a finite extension of $\mathbb Q$), $s>0$ be the number of its real embeddings and $2t>0$ be the number of its complex embeddings. One can easily prove that  for each $s$ and $t$ there exists a field $K$ which has these numbers of real and complex embeddings (see e.g. \cite{_Oeljeklaus_Toma_}).

\definition The \textit{ring of algebraic integers} $O_K$ is a subring of $K$ that consists of all roots of polynomials with integer coefficients which lie in $K$.
\textit{Unit group} $O^*_K$ is the multiplicative subgroup of invertible elements of $O_K$.

\hfill

Let $m$ be $s+t$. Let $\sigma_1, \ldots, \sigma_s$ be real embeddings of the field $K$, $\sigma_{s+1}, \ldots, \sigma_{s+2t}$ be complex embeddings
such that $\sigma_{s+i}$ and $\sigma_{s+t+i}$ are complex conjugate for each $i$ from $1$ to $t$. Now we can define a map 
$l: O_K^*\rightarrow \mathbb R^m$ where $l(u)=(\ln|\sigma_1(u)|, \ldots, \ln|\sigma_s(u)|, 2\ln|\sigma_{s+1}(u)|,
\ldots, 2\ln|\sigma_m(u)|)$.

Denote $O_K^{*, +}=\{a\in O_K^*:\sigma_i(a)>0, i=1, \ldots, s\}$. Let us consider following definitions:

\hfill

\definition A \textit{lattice} $\Lambda$ in $\mathbb R^n$ is a discrete additive subgroup such that $\Lambda\otimes\mathbb R=\mathbb R^n$.

\definition \cite{_Oeljeklaus_Toma_} The group $U\subset O_K^{*, +}$ of rank $s$ is called \textit{admissible for the field} $K$ if the projection of $l(U)$ to the first $s$ components is a lattice in $\mathbb R^s$.

Consider a linear space $L=\{x\in\mathbb R^m \mid \sum_{i=1}^m x_i=0\}$. The projection of $L\subset R^m$ to the first $s$ coordinates is surjective, because $s<m$.
Using the Dirichlet unit theorem (see e.g. \cite{_Milne09_}) one can prove that $l(O_K^{*, +})$ is a full lattice in $L$. Therefore there exists a group $U$ that is admissible.

\hfill

Let $\mathbb H=\{z\in\mathbb C \mid \im z>0\}$. Let $U\subset O_K^{*, +}$ be a group which is admissible for $K$. The group $U$ acts on $O_K$ multiplicatively. This defines a structure of 
semidirect product $U':=U\ltimes O_K$. Define the action of $U'$ on $\mathbb H^s\times\mathbb C^t$ as follows. The element $u\in U$ acts on $\mathbb H^s\times \mathbb C^t$ mapping $(z_1, \ldots, z_m)$ to $(\sigma_1(u)z_1, \ldots, \sigma_m(u)z_m)$. Since 
$U$ lies in $O_K^{*, +}$, the action $U$ on the first $s$ coordinates preserves $\mathbb H$.

The additive group $O_K$ acts on $\mathbb H^s\times \mathbb C^t$ by parallel translations: $a\in O_K$ is mapping $(z_1, \ldots, z_m)$ to 
$(\sigma_1(a)+z_1, \ldots, \sigma_m(a)+z_m)$. Since the first $s$ embeddings are real, this action preserves $\mathbb H$ in the first $s$ coordinates.

One can see that $(u, a)\in U\ltimes O_K$ maps $(z_1, \ldots, z_m)$ to $(\sigma_1(u)z_1+\sigma_1(a), \ldots, \sigma_m(u)z_m+\sigma_m(a))$. One can easily show that this action
is compatible with the group operation in the semidirect product.

\hfill

\definition An \textit{Oeljeklaus-Toma manifold} is the quotient of $\mathbb H^s\times\mathbb C^t$ by the action of the group $U\ltimes O_K$, which was defined above.

This quotient exists because $U\ltimes O_K$ acts properly discontiniously on $\mathbb H^s\times\mathbb C^t$. Additionally $\mathbb H^s\times\mathbb C^t/U\ltimes O_K$ is a compact complex manifold. To prove it, let $U$ be admissible for $K$. The quotient $\mathbb H^s\times\mathbb C^t/O_K$ is obviously diffeomorphic to the trivial toric bundle $(\mathbb R_{>0})^s\times (S^1)^n$. The group $U$ acts properly discontinuously on the base $(\mathbb R_{>0})^s$. Therefore it acts properly discontinuously on $\mathbb H^s\times\mathbb C^t/O_K$. Also, the groups $U$ and $O_K$ act holomorphically on $\mathbb H^s\times\mathbb C^t$. Therefore the quotient has a holomorphic structure.

\section{Curves on the Oeljeklaus-Toma manifolds}
In this section we shall prove that there are no complex curves on the Oeljeklaus-Toma manifolds, just as on
Inoue surfaces of type $S_M$ (see \cite{_Inoue_}).

\subsection{The exact semipositive (1,1)-form on the Oeljeklaus-Toma manifold}

The $(1,1)$-form we will be using was previously introduced in the paper \textit{Subvarieties in Oeljeklaus-Toma 
manifolds} by Ornea and Verbitsky. Authors use this form to prove that Oeljeklaus-Toma manifolds with $t=1$
(that means that the corresponding number field has only two complex embeddings) do not contain any 
submanifolds.
Our result works for all Oeljeklaus-Toma manifolds, but we
consider only curves instead of submanifolds of any dimension. Later on we will explain how our method differs
from the one in \cite{_Ornea_Verbitsky_} and how this implies the difference between our results.

Define the notion of $(1,1)$-form. 
Let $M$ be a smooth complex manifold, $z_1, \ldots, z_n$ --- local complex coordinates in the open neighborhood 
of the point $y\in M$.

\definition A \textit{$(1,1)$-form} on a complex manifold $M$ is a 2-form $\omega$, such that $\omega(Iu,v)=
-\omega(u,Iv)=\sqrt{-1}\omega(u,v)$ for each $u,v\in T_yM$, where $I$ is the almost complex structure on $M$.

\definition A $(1, 1)$-form $\omega$ on a complex manifold $M$ is \textit{semipositive} 
if $\omega(u, Iu)\geqslant 0$ for each tangent vector $u\in T_yM$.

\hfill

As in \cite{_Ornea_Verbitsky_}, we consider a certain semipositive $(1, 1)$-form on the Oeljeklaus-Toma manifold
$M=\mathbb H^s\times\mathbb C^t/(U\ltimes O_K)$. We introduce a $(1, 1)$-form $\tilde\omega$ on $\tilde M=\mathbb H^s\times\mathbb C^t$ which is preserved by the action of the group $\Gamma=(U\ltimes O_K)$ and since then it would be a $(1,1)$-form on $M$.

Let $(z_1, \ldots, z_m)$ be complex coordinates on $\tilde M$. Define $\varphi(z)=\Pi_{i=1}^s \im(z_i)^{-1}$.
Since the first $s$ components of $\tilde M$ correspond to upper half-planes $\mathbb H\subset\mathbb C$, this function is positive on $\tilde M$. 

\hfill

Let us now consider the form $\tilde\omega=\sqrt{-1}\partial\bar\partial\log\varphi$. Using standard coordinates on $\tilde M$ one can write this form as $\tilde\omega=\sqrt{-1}\sum_{i=1}^s\frac{dz_i\wedge d\bar z_i}{4(\im z_i)^2}$. Therefore $\tilde\omega$ is a semipositive $(1,1)$-form on $\tilde M$.

Let us show that this form is $\Gamma$-invariant.

\hfill

The group $\Gamma$ is a semidirect product of the additive group $O_K$ and the multiplicative group $U$.
The additive group acts on the first $s$ components of $\tilde M$ (which correspond to upper half-planes $\mathbb H\subset\mathbb C$) by translations along the real line. Therefore it does not change $\im z_i$ for $i=1\ldots s$. Hence the function $\log\varphi$ is preserved by the action of the additive component.

The multiplicative component acts on the first $s$ coordinates of $\tilde M$ by multiplying them by a real number (since the first $s$ embeddings of the number field $K$ are real). Then every $\im z_i$ is multiplied by a real number and so there is a real number added to $\log (\im z_i)$. Since $\log\varphi(z)=-\sum_{i=1}^s\log (\im z_i)$, there is a real number added to $\log\varphi$. The operator $\bar\partial$ is zero on the constants, so $\tilde\omega=\sqrt{-1}\partial\bar\partial\log\varphi$ is preserved by action of the group $\Gamma$.

Since the $(1,1)$-form $\tilde\omega$ is $\Gamma$-invariant it is the pullback of $(1,1)$-form $\omega$ on the Oeljeklaus-Toma manifold $M=\tilde M/\Gamma$.

Let us now show that the form $\tilde\omega$ is exact on $\tilde M$. For that we define the operator $d^c$.

\definition Define the \textit{twisted differential} $d^c=I^{-1}dI$ where $d$ is a De Rham differential and $I$ is the almost complex structure.

Since $dd^c=2\sqrt{-1}\partial\bar\partial$ (see \cite{_Griffits-Harris_}), one can see that $\tilde\omega=\sqrt{-1}\partial\bar\partial\log\varphi=\frac{1}{2}dd^c\log\varphi$ and so $\tilde\omega$ is exact as a form on $\tilde M$. Also since the operator $d^c$ vanishes on constants the form $d^c\log\varphi$ is $\Gamma$-invariant, so $\omega$ is exact on $M$.

\subsection{The $(1, 1)$-form $\omega$ and curves on the Oeljeklaus-Toma manifold}

Since the form $\omega$ on the manifold $M$ is semipositive, its integral on any complex curve $C\subset M$ is nonnegative. The form $\omega$ is exact. Hence Stokes' theorem implies that its integral on any complex curve vanishes. Therefore if $C\subset M$ is a closed complex curve,  $\omega$ vanishes on it.

To find out on which curves $\omega$ vanishes, let us define the zero foliation of the form $\omega$.

\hfill

\definition An \textit{involutive distribution (or foliation)} on $M$ is a subbundle $B\subset TM$ of the tangent bundle that is closed under the Lie bracket: $[B, B]\subset B$.

\definition A \textit{leaf of a foliation} $B$ is a connected submanifold of $M$ such that its dimension is equal to $\dim B$ and that is tangent to $B$ at every point.

\theorem (Frobenius) Let $B\subset TM$ be an involutive distribution. Then for each point of the manifold $M$, there is exactly one leaf of this distribution that contains this point
(see e.g. \cite{_Boothby_} Section IV. 8. Frobenius Theorem).

\theorem Let $N\subset M$ be a connected submanifold such that its tangent space at every point lies in a foliation $F\subset TM$. Then $N$ lies in a leaf of the foliation $F$ (see e.g. \cite{_Boothby_} Section IV. 8. Theorem 8.5).

\hfill

\definition The \textit{zero foliation} of a semipositive $(1,1)$-form $\omega$ on $M$ is the subundle of $TM$ that consists of tangent vectors $u\in T_yM$ such that $\omega(u, Iu)=0$, where $I$ is the almost complex structure on $M$.

\hfill

Consider the zero foliation of $\tilde\omega$ on $\tilde M$.

The form $\tilde\omega$ is strictly positive on each vector $v=(z_1,\ldots,z_m)$ such that at least one of $z_i$ for $i=1,\ldots, s$ is nonzero. Such a vector cannot be tangent to a leaf of the zero foliation. Therefore on each leaf of the zero foliation of the form $\tilde\omega$ the first $s$ coordinates are constant.

Hence a leaf of the zero foliation of $\tilde\omega$ on $\tilde M$ is isomorphic to $\mathbb C^t$.

\hfill

Let us now consider the zero foliation of $\omega$ on $M$.

We show that the non-trivial image of the action of $\Gamma$ on any leaf $L$ of the zero foliation of the form $\tilde\omega$ does not intersect with $L$.

One can see that $L$ is $(z_1, \ldots, z_s)\times\mathbb C^t$ for some fixed $(z_1, \ldots, z_s)$. Therefore, for any $\gamma\in\Gamma$ such that $L\cap \gamma(L)\neq\emptyset$, the first $s$ coordinates of the points in $L$ coincide with the first $s$ coordinates of the points in $\gamma(L)$. Then for such $\gamma$ have the following system of equations:

$$\sigma_i(u)z_i+\sigma_i(a)=z_i,\quad i=1\ldots s,$$
where $\gamma=(u, a)$.

These equations imply that $z_i=\frac{\sigma_i(a)}{1-\sigma_i(u)}$. Therefore $z_i$ are real but $\mathbb H$ does not have real elements.

We showed that $L\cap \gamma(L)=\emptyset$ for every $\gamma\neq 1$ in $\Gamma$.

Since $\omega$ vanishes on each compact curve $C\subset M$, each curve is contained in some leaf of the zero foliation of $\omega$. Since $\tilde\omega$ is $\Gamma$-invariant, each leaf of the zero foliation of $\omega$ on $M$ is isomorphic to a component of the leaf of the zero foliation of $\tilde\omega$ on $\tilde M$. Therefore, it is isomorphic to $\mathbb C^t$. And $\mathbb C^t$ does not contain any compact complex submanifolds.

We proved the following theorem:

\hfill

\theorem There are no compact complex curves on the Oeljeklaus-Toma manifolds.

\section{Closing remarks}

Let us now briefly explain the connection between our work and \cite{_Ornea_Verbitsky_}. As in \cite{_Ornea_Verbitsky_} we use the zero foliation of a certain $(1,1)$-form. The leaves of this zero foliation are $t$-dimensional complex manifolds. In \cite{_Ornea_Verbitsky_} authors consider $t=1$, and submanifolds of any dimension. We consider any $t\geq 1$, and submanifolds of dimension 1. In both works one uses the semipositivity of the $(1,1)$-form to prove that a submanifold and a leaf of the zero foliation, which contains a point in this submanifold, could not intersect transversely. In both cases, one of the manifolds is of dimension one. Non-transversality implies that it is contained in the second one. Authors of \cite{_Ornea_Verbitsky_} use the fact that the leaves of the zero foliation are Zariski dense, and therefore could not be contained in any submanifold. In our work we prove that each leaf is isomorphic to $\mathbb C^t$ and therefore could not contain any submanifolds.

The arguments used in our work and \cite{_Ornea_Verbitsky_} could not be applied to the case of higher dimensions.

It seems obvious that there should be Oeljeklaus-Toma manifols which contain submanifolds --- other Oeljeklaus-Toma manifolds, corresponding to smaller number fields, but we do not have a formal proof.

In \cite{_Oeljeklaus_Toma_} it was proved that Oeljeklaus-Toma manifolds do not admit non-trivial meromorphic functions. Therefore all the divisors on these manifolds are fixed in their linear systems (which are 0-dimensional). It is conjectured that Oeljeklaus-Toma manifolds admit no divisors.

\hfill

\hfill

\begin{small}
\noindent{\sc Sima Verbitsky\\
{\sc Moscow State University, Faculty of Mathematics and Mechanics\\
GSP-1 1, Leninskie gory, 119991 Moscow, Russia.}\\
\tt sverb57@gmail.com}
\end{small}


{\small
\begin{thebibliography}{XXXI}

\bibitem[Boo]{_Boothby_}
Boothby W.M. \textit{An Introduction to Differentiable Manifolds and Riemannian Geometry}. Academic Press, San Diego, California, 2003.

\bibitem[G--H]{_Griffits-Harris_}
Griffiths Ph., Harris J. \textit{Principles of Algebraic Geometry}. Wiley-Interscience, 1994.

\bibitem[I]{_Inoue_}
Inoue M. \textit{On surfaces of Class} VII$_0$, Invent. Math. 24 (1974), 269-310.

\bibitem[Mil08]{_Milne08_}
Milne J.S. \textit{Fields and Galois Theory, September 2008}.\\
This paper can be found on http://www.jmilne.org/math/CourseNotes/ft.html, version 4.21

\bibitem[Mil09]{_Milne09_}
Milne J.S. \textit{Algebraic Number Theory, April 2009}.\\
This paper can be found on http://www.jmilne.org/math/CourseNotes/ant.html, version 3.02

\bibitem[O--T]{_Oeljeklaus_Toma_}
Oeljeklaus K., Toma M. \textit{Non-K\" ahler compact complex manifolds associated to number fields}. Ann. Inst. Fourier 55 (2005), 1291-1300.

\bibitem[O--V]{_Ornea_Verbitsky_}
Ornea L., Verbitsky M. \textit{Subvarieties in Oeljeklaus-Toma manifolds}.

\bibitem[P--V]{_Parton_Vuletescu_}
Parton M., Vuletescu V. \textit{Examples of non-trivial rank in locally conformal K\" ahler geometry}. Math. Z. (2010), DOI 10.1007/s00209-010-0791-5, arXiv:1001.4891.

\bibitem[R]{_Raghunathan_}
Raghunathan M.S. \textit{Discrete subgroups of Lie groups}. Springer 1972.

\bibitem[V]{_Voisin_}
Voisin C. \textit{Hodge Theory and Complex Algebraic Geometry Volume 1}.
Cambridge University Press, 2002.
\end{thebibliography}
}
\end{document}